\documentclass[a4paper,11pt]{amsart}

\usepackage[left=2cm,top=3cm,bottom=2cm,right=2cm]{geometry}

\usepackage{amsthm}
\usepackage{amssymb}
\usepackage{amsmath}
\usepackage{mathtools}
\usepackage{pdfpages}
\usepackage{changes}
\usepackage{exscale,cite,color,amsopn}
\usepackage[colorlinks=true,pdftex,unicode=true,linktocpage,bookmarksopen,hypertexnames=true]{hyperref}
\usepackage{aliascnt}
\usepackage{colortbl}
\usepackage{enumitem}
\usepackage{verbatim}
\usepackage{tikz-cd}
\usepackage[capitalise]{cleveref}

\newcounter{namedcondition}
\crefname{namedcondition}{condition}{conditions}
\Crefname{namedcondition}{Condition}{Conditions}

\newcommand{\nameditem}[2]{
    \item[#2]
    \begingroup
        \renewcommand{\thenamedcondition}{#2}
        \refstepcounter{namedcondition}
        \label{#1}
    \endgroup
}

\usetikzlibrary{trees}

\renewcommand{\leq}{\leqslant}
\renewcommand{\geq}{\geqslant}

\DeclareMathOperator{\id}{id}

\DeclareMathOperator{\Ret}{Ret}

\DeclareMathOperator{\Sym}{Sym}

\DeclareMathOperator{\ret}{ret}

\newcommand{\N}{\mathbb{N}}

\newcommand{\F}{\mathbb{F}}
\newcommand{\Q}{\mathbb{Q}}

\crefname{table}{table}{tables}

\newcommand{\genrel}[2]{\left\langle \ #1 \ \vline \ #2 \ \right\rangle}

\makeatletter
\numberwithin{equation}{section}
\numberwithin{figure}{section}
\numberwithin{table}{section}
\newtheorem{thm}{Theorem}[section]
\newtheorem*{thm*}{Theorem}
\newtheorem{lem}[thm]{Lemma}

\newtheorem{pro}[thm]{Proposition}

\newtheorem*{convention*}{Convention}

\theoremstyle{definition}
\newtheorem{defn}[thm]{Definition}
\newtheorem{problem}[thm]{Problem}
\newtheorem*{problem*}{Problem}
\newtheorem*{conjecture*}{Conjecture}

\makeatother

\title[Structure groups without the UPP]{Structure groups of solutions to the Yang--Baxter equation without the unique product property}
\author{Carsten Dietzel}
\date{\today}

\address[Carsten Dietzel]{Department of Mathematics and Data Science, Vrije Universiteit Brussel, Pleinlaan 2, 1050 Brussel, Belgium}
\email{Carsten.Dietzel@vub.be}

\begin{document}

\begin{abstract}
    The unique product property is a group-theoretic property which provides an important combinatorial obstruction to the existence of special elements in group rings of torsionfree groups. However, it is surprisingly difficult to construct concrete examples of torsionfree groups without the unique product property.
    Due to their abundance of quadratic relations, structure groups of set-theoretic solutions to the Yang--Baxter equation are natural candidates for torsionfree groups without the unique product property. Acri, Lutowski and Vendramin found such structure groups using algorithms that detect subgroups isomorphic to the Promislow group, a prototypical example of a torsionfree group without the unique product property.
    Combining string rewriting techniques with the Garside structure of structure groups, we prove the existence of set-theoretic solutions to the Yang--Baxter equation whose structure group neither satisfies the unique product property nor contains the Promislow group. 
\end{abstract}

\maketitle

\section*{Introduction}

A group $G$ is said to satisfy the \emph{unique product property} (\emph{UPP}, for short) if for any two finite, nonempty subsets $A,B \subseteq G$, there is an element $g \in G$ that can be represented as $g = ab$ with $a \in A$, $b \in B$ in exactly one way. Note that if $G$ has the UPP then $G$ is torsionfree.

Groups with the UPP were first investigated by Rudin and Schneider \cite{Rudin_Schneider_UPP} under the name \emph{$\Omega$-groups} in the context of idempotents in group rings, addressing a problem of Kaplansky \cite{KaplanskyProblems}. We now recall Higman's \emph{unit conjecture} \cite{HigmanThesis,HigmanUnit}, often attributed to Kaplansky due to being mentioned in a problem list \cite{KaplanskyRevisited}.

\begin{conjecture*}
    Let $k$ be a field, and let $G$ be a torsionfree group. Then each unit $u$ in the group ring $k[G]$ is of the form $u = rg$ with $r \in k^{\times}$, $g \in G$.
\end{conjecture*}

It is easy to see that in a potential counterexample, $G$ cannot be a unique product group. However, finding torsionfree groups without the UPP is a more delicate problem. The first examples of such groups were constructed by Rips and Segev \cite{RS_torsionfree_UPP}. Promislow \cite{Promislow_UPP} subsequently showed that there is a more natural example. The \emph{Promislow group} (also called a \emph{Hantzsche--Wendt group} in the theory of space groups \cite{HW_Raumformen}) is a three-dimensional crystallographic group that has a presentation by generators and relations given by
\[
P = \genrel{x,y}{x^2yx^2 = y;\ y^2xy^2 = x}_{\mathrm{gr}},
\]
where the subscript $\mathrm{gr}$ indicates that the presentation is to be understood in terms of groups. 

In 2021, Gardam \cite{Gardam_Kaplansky_F2} disproved the unit conjecture over $\mathbb{F}_2$, the field with two elements, by providing explicit nontrivial units in the group ring $\F_2[P]$. Shortly thereafter, Murray \cite{Murray_Kaplansky} showed that the group rings $\F_p[P]$ also contain nontrivial units for all primes $p$. Finally, in 2024, Gardam \cite{Gardam_Kaplansky_C} proved that $\mathbb{C}[P]$ also contains nontrivial units, thus settling the unit conjecture in characteristic $0$. To the author's knowledge, it is an open question if the unit conjecture holds over the prime field $\Q$.

In order to understand the failure of the UPP, it is natural to consider the structure groups of \emph{involutive, nondegenerate set-theoretic solutions to the Yang--Baxter equation} (or \emph{solutions}, for short). A solution $(X,r)$ is given by a finite set $X$, together with an involution $r: X^2 \to X^2$ that satisfies the following equation in $\Sym_{X^3}$:
\begin{equation}
    (r \times \id_X)(\id_X \times r)(r \times \id_X) = (\id_X \times r)(r \times \id_X)(\id_X \times r)
\end{equation}
(see \cref{subsec:the_YBE} for a precise definition of the nondegeneracy condition). An important group-theoretic invariant of solutions is their \emph{structure group}, which is presented, in terms of generators and relations, as
\[
G(X,r) = \genrel{X}{wx = yz \textnormal{ whenever } r(w,x) = (y,z)}_{\mathrm{gr}}.
\]
Since $G(X,r)$ is torsionfree while having an abundance of nontrivial quadratic relations, it is natural to pose the following question.

\begin{problem*}
    Let $(X,r)$ be a solution to the Yang--Baxter equation. When does $G(X,r)$ satisfy the UPP?
\end{problem*}

It has been shown by Bachiller, Ced\'o and Vendramin \cite{BCV_Left_Orderable} that $G(X,r)$ is left-orderable if and only if $(X,r)$ is of \emph{finite multipermutation level} (see \cref{subsec:the_YBE} for a precise definition). As left-orderability implies the UPP\footnote{But not every group satisfying the UPP is left-orderable, see \cite[Appendix B]{RK_diffuse_groups}, for example.}
, $G(X,r)$ can only violate the UPP if $(X,r)$ is not of finite multipermutation level.

Due to results of Castelli, Catino and Pinto \cite{CCP_indecomposable}, it is known that for \emph{certain} surmorphisms $f: (X,r) \twoheadrightarrow (Y,s)$ of solutions, the UPP for $G(X,r)$ implies the UPP for $G(Y,s)$. The conditions imposed on $f$ ensure that $G(Y,s)$ embeds into $G(X,r)$.

In this article, we show that these additional conditions are unnecessary.

\begin{thm*}
    Let $(X,r)$ and $(Y,s)$ be solutions to the Yang--Baxter equation and suppose that there is a surmorphism $f: (X,r) \twoheadrightarrow (Y,s)$. If $G(X,r)$ satisfies the UPP, then $G(Y,s)$ satisfies the UPP.
\end{thm*}

Computational results of Acri, Lutowski and Vendramin \cite{ALV_UPP} show that for all but $8$ isomorphism classes of solutions $(X,r)$ with $|X| \leq 8$, the group $G(X,r)$ satisfies the UPP if and only if $(X,r)$ is of finite multipermutation level. By an application of our theorem, we are able to extend this result to all solutions of size at most $8$.

\begin{thm*}
    Let $(X,r)$ be a solution to the Yang--Baxter equation with $|X| \leq 8$. Then $G(X,r)$ satisfies the UPP if and only if $(X,r)$ is of finite multipermutation level.
\end{thm*}

We remark here that the authors of \cite{ALV_UPP} identified those solutions $(X,r)$ to the Yang--Baxter equation with $|X| \leq 8$ such that $G(X,r)$ contains the Promislow group as a subgroup, and therefore fails to satisfy the UPP. Amongst those solutions that are not of finite multipermutation level, there are $8$ exceptional isomorphism classes for which $G(X,r)$ does not contain the Promislow group as a subgroup.

As a byproduct, our proof yields $8$ explicit examples of solutions $(X,r)$ for which $G(X,r)$ does not satisfy the UPP, while also not containing the Promislow group.

\subsection*{Acknowledgements and Funding}

The author of this article expresses his gratitude to FWO Vlaanderen which supports him by means of a Senior Postdoctoral Fellowship, grant number 1244625N. This work was partially supported by the project OZR3762 of Vrije Universiteit Brussel.

\subsection*{AI disclosure statement}

During the preparation of the manuscript, ChatGPT-6.2 (OpenAI) has exclusively been used for proofreading the manuscript in the final stage. The mathematical ideas and arguments presented here have been developed without AI assistance. 

\section{Preliminaries} \label{sec:preliminaries}

\subsection{Monoids} \label{subsec:monoids}

We assume that the reader is familiar with the theory of monoids and groups defined by generators and relations. However, since we will make heavy use of string rewriting in the context of monoids, we introduce the necessary notation by recapitulating the basics of monoids defined by generators and relations. Further background on monoid and group presentations can be found in \cite{Higgins_Semigroups} and \cite{LS_Combinatorial_Group_Theory}, respectively.

Recall that for a set $X$, the \emph{free monoid} over $X$ is defined as the set $X^{\ast} = \{\varepsilon \} \cup \bigcup_{k = 1}^{\infty} X^k$. An element $(x_1,x_2,\ldots,x_k) \in X^{\ast}$ is denoted as a string $x_1x_2\ldots x_k$, $\varepsilon$ representing the empty string, and the monoid operation on $X^{\ast}$ is given by concatenation of strings. Accordingly, we refer in this context to the elements of $X$ as \emph{letters} and to the elements of $X^{\ast}$ as \emph{words}.

It is well-known that $X^{\ast}$ satisfies the following universal property. Let $\iota_X: X \hookrightarrow X^{\ast}$; $x \mapsto x$ be the canonical embedding. If $M$ is a monoid, and $\varphi: X \to M$ is any mapping of sets, then there is a unique homomorphism of monoids $\varphi^{\ast}: X^{\ast} \to M$ with $\varphi^{\ast} \iota_X = \varphi$. This homomorphism is explicitly given by $\varphi^{\ast}(x_1x_2\ldots x_k) = \varphi(x_1)\varphi(x_2)\ldots \varphi(x_k)$.

A set of \emph{relations} in $X^{\ast}$ is a subset $R \subseteq X^{\ast} \times X^{\ast}$ that we will always assume to be \emph{symmetric}, that is $(v,w) \in R$ implies $(w,v) \in R$ for all $v,w \in X^{\ast}$. However, we do not necessarily assume that it is also transitive.

Given a set of relations $R \subseteq X^{\ast} \times X^{\ast}$, we write $v \sim_R w$ for $v,w \in X^{\ast}$ if there are factorizations $v = v_1 v_2 v_3$ and $w = w_1 w_2 w_3$, with possibly empty factors, such that $v_1 = w_1$, $v_3 = w_3$ and $(v_2,w_2) \in R$. Furthermore, we write $v \approx_R w$ if for some $n\geq 0$, there is a sequence
\[
v = v_0 \sim_R v_1 \sim_R \ldots \sim_R v_n = w
\]
with all $v_i \in X^{\ast}$ ($0 \leq i \leq n$). It is easy to see that $\approx_R$ is a congruence relation on $X^{\ast}$, and we define $\genrel{X}{R}_{\mathrm{mon}} = X^{\ast}/\!\approx_R$.  In fact, $\approx_R$ is the smallest congruence relation on $X^{\ast}$ such that $R \subseteq \approx_R$. As usual, we suggestively write a relation $(v,w) \in R$ as $v = w$, if it appears in a presentation of the form $\genrel{X}{R}_{\mathrm{mon}}$.

Given a map $f: X \to Y$, we denote by $f^{\ast}: X^{\ast} \to Y^{\ast}$ the canonical monoid homomorphism given by $f^{\ast}(\varepsilon) = \varepsilon$ and $f^{\ast}(x_1x_2\ldots x_k) = f(x_1) f(x_2) \ldots f(x_k)$. Obviously not every monoid homomorphism from $X^{\ast}$ to $Y^{\ast}$ is of this form.

In this work, we will need the following simple lemma.

\begin{lem} \label{lem:respect_relations}
    Let $M = \genrel{X}{R}_{\mathrm{mon}}$ and $N = \genrel{Y}{S}_{\mathrm{mon}}$ and let $f: X \to Y$ be a map such that $(f^{\ast}(v),f^{\ast}(w)) \in S$ whenever $(v,w) \in R$. Then $f^{\ast}$ induces a well-defined homomorphism $\overline{f}: M \to N$; $\overline{v} \mapsto \overline{f^{\ast}(v)}$.
\end{lem}

\begin{proof}
    Let $F:X^{\ast} \to N$; $v \mapsto \overline{f^{\ast}(v)}$ and consider the congruence relation $v \sim w \Leftrightarrow F(v) = F(w)$. By assumption $F(v) = F(w)$ for any $(v,w) \in R$, and as $\approx_R$ is the smallest congruence containing $R$, the homomorphism $F$ factors through $X^{\ast}/\!\!\approx_R = M$. 
\end{proof}

\subsection{Garside groups} \label{subsec:garside_groups}

In the following, we assume familiarity with the concept of a lattice. The classical definition of a Garside group is rather technical, so we use a more convenient definition in terms of right $\ell$-groups, which has been shown by Rump to be equivalent to the classical definition (see \cite[Sections 1+6]{Rump-Geometric-Garside}).

Let $G$ be a group that is equipped with a lattice structure whose join and meet operations are denoted by $\vee$ and $\wedge$, respectively. Then $G = (G,\vee,\wedge)$ is called a \emph{right $\ell$-group} if for all $x,y,z \in G$, we have
\[
(x \vee y)z = xz \vee yz \ ; \ (x \wedge y)z = xz \wedge yz.
\]
Given a right $\ell$-group $G$, its \emph{positive monoid} is defined as the submonoid $G^+ = \{g \in G: g \geq 1 \}$. Similarly, the \emph{negative monoid} is defined as the submonoid $G^- = \{g \in G: g \leq 1 \}$.

A right $\ell$-group $G$ is \emph{noetherian} if all ascending sequences $g_1 \leq g_2 \leq \ldots$ with $g_k \in G^-$ for all $k \geq 1$ become stationary, and similarly, if all descending sequences $g_1 \geq g_2 \geq \ldots$ with $g_k \in G^+$ for all $k \geq 1$ become stationary.

Furthermore, an element $\Delta \in G^+$ is called a \emph{strong order unit} if $\Delta$ is \emph{normal}, meaning that $\Delta(x \vee y) = \Delta x \vee \Delta y$ and $\Delta(x \wedge y) = \Delta x \wedge \Delta y$ hold for all $x,y \in G$, and if for all $g \in G$, there exists an integer $k$ such that $g \leq \Delta^k$.

We can now state one possible definition of a \emph{Garside group}.

\begin{defn}
    A \emph{Garside group} is a noetherian right $\ell$-group with a distinguished strong order unit $\Delta \in G^+$, the \emph{Garside element}, such that the \emph{Garside interval} $[1,\Delta] = \{ g \in G^+ : g \leq \Delta \}$ is finite.

    The positive monoid $G^+$ is called the \emph{Garside monoid} of $G$.
\end{defn}

We will make use of the following property of Garside groups.

\begin{lem} \label{lem:shifting_subsets_into_G+}
    Let $G$ be a Garside group with Garside element $\Delta$. Then for all finite subsets $S \subseteq G$, there is an integer $k$ such that $\Delta^k S \subseteq G^+$. Similarly, there is an integer $l$ such that $S \Delta^l \subseteq G^+$.
\end{lem}

\begin{proof}
    If $S = \emptyset$, the statement is trivial, so suppose that $S \neq \emptyset$. Let $g = \bigwedge S$. Due to the normality of $\Delta$, it is sufficient for both statements to show that there is an integer $k$ with $g \geq \Delta^{-k}$. Due to the normality of $\Delta$, this is equivalent to $\Delta^k \geq g^{-1}$, and the strong order unit property of $\Delta$ implies the existence of such an integer $k$.
\end{proof}

\subsection{The Yang--Baxter equation} \label{subsec:the_YBE}

In this subsection, we recall the basic facts and definitions concerning set-theoretic solutions to the Yang--Baxter equation. All relevant material, except for \cite{Chouraqui_Garside}, can be found in the seminal paper \cite{ESS_YangBaxter}.

Let $X$ be a set and let $r: X^2 \to X^2$; $(x,y) \mapsto r(x,y) = (\sigma_x(y),\tau_y(x))$ be a map. One says that $r$ satisfies the \emph{set-theoretic Yang--Baxter equation} if $r$ is bijective and if the following identity of maps $X^3 \to X^3$ is satisfied:
\begin{equation} \label{eq:YBE}
    (r \times \id_X)(\id_X \times r)(r \times \id_X) = (\id_X \times r)(r \times \id_X)(\id_X \times r). \tag{YBE}
\end{equation}
If the coordinate maps $\sigma_x$, $\tau_x$ are bijective for all $x \in X$, then $r$ is said to be \emph{nondegenerate}. Furthermore, $r$ is called \emph{involutive} if $r^2 = \id_{X^2}$.

\begin{convention*}
    In this article, when speaking of a \emph{solution to the Yang--Baxter equation} (or \emph{solution}, for short), we will mean a nonempty, finite, involutive and nondegenerate solution to \eqref{eq:YBE}.
\end{convention*}

Given solutions $(X,r)$ and $(Y,s)$, a map $f: X \to Y$ is called a \emph{morphism} of solutions if $s(f \times f) = (f \times f)r$. If $f: (X,r) \to (Y,s)$ is a surjective morphism of solutions, we will briefly call it a \emph{surmorphism}. Furthermore, given a solution $(X,r)$, we call a nonempty subset $Y \subseteq X$ a \emph{subsolution} if $r(Y^2) = Y^2$. In such a case, $(Y,r) = (Y,r|_{Y^2})$ is again a solution to the Yang--Baxter equation.

A particularly important example of a surmorphism is the \emph{retraction map} of a solution. Given a solution $(X,r)$ with $r(x,y) = (\sigma_x(y),\tau_y(x))$, the \emph{retraction relation} is defined by $x \sim y \Leftrightarrow \sigma_x = \sigma_y$. It can be shown that the set $\overline{X} = X/\!\!\sim$ with the induced map $\overline{r}: \overline{X}^2 \to \overline{X}^2$; $([x],[y]) \mapsto ([\sigma_x(y)],[\tau_y(x)])$ defines a solution in a well-defined manner. This solution is called the \emph{retraction} of $(X,r)$ and is denoted by $\Ret(X,r) = (\overline{X},\overline{r})$. Furthermore, the \emph{retraction map} $\ret: (X,r) \twoheadrightarrow \Ret(X,r)$; $x \mapsto [x]$ is a surmorphism of solutions. If the retraction map of $(X,r)$ is a bijection, we call $(X,r)$ \emph{irretractable}, and \emph{retractable} otherwise.

For a solution $(X,r)$, one defines inductively $\Ret^0(X,r) = (X,r)$ and $\Ret^{k+1}(X,r) = \Ret(\Ret^k(X,r))$ ($k \geq 0$), and if there exists a $k \geq 0$ such that $\Ret^k(X,r)$ is a solution on a set with one element, then $(X,r)$ is said to be of \emph{finite multipermutation level}.

We now recall the definition of the structure monoid and the structure group of a solution.

\begin{defn}
    Let $(X,r)$ be a solution to the Yang--Baxter equation. Its \emph{structure monoid} is defined as
    \begin{equation}
        M(X,r) = \genrel{X}{xy = wz \textnormal{ whenever } r(x,y) = (w,z)}_{\mathrm{mon}}
    \end{equation}
    Similarly, its \emph{structure group} is given by
    \begin{equation}
        G(X,r) = \genrel{X}{xy = wz \textnormal{ whenever } r(x,y) = (w,z)}_{\mathrm{gr}}
    \end{equation}
\end{defn}

A remarkable property of the structure group is that is has a Garside structure. In \cref{sec:structure_groups_and_the_upp}, this will allow us to translate statements about $G(X,r)$ to statements about $M(X,r)$.

\begin{thm} \label{thm:chouraqui_garside}
    Let $(X,r)$ be a solution to the Yang--Baxter equation. Then the canonical homomorphism $\iota: M(X,r) \to G(X,r)$ identifies $M(X,r)$ with the Garside monoid of a Garside structure on $G(X,r)$. In particular, $\iota$ is injective.
\end{thm}

\begin{proof}
    \cite[Theorem 3.3]{Chouraqui_Garside}.
\end{proof}

\section{Essential maps for monoid presentations} \label{sec:essential_maps}

In this section, we investigate the UPP in the more general setting of monoids.

If $M$ is a monoid, and $A,B \subseteq M$, we write $AB = \{ ab : a \in A, b \in B \} \subseteq M$ for the \emph{product set} of $A$ and $B$ in $M$.

\begin{defn}
    Let $M$ be a monoid. We say that $M$ satisfies the \emph{UPP} if for any finite subsets $\emptyset \neq A,B \subseteq M$, there is an element $x \in AB$ that is representable in exactly one way as $x = ab$ with $a \in A$, $b \in B$. 
\end{defn}

In this section, our aim is to show that under suitable conditions, the UPP descends along surmorphisms of monoids $M \twoheadrightarrow N$, given that the homomorphism is induced by an \emph{essential} map with respect to the defining relations. The reader familiar with algebraic topology will notice the analogy between essential maps and fibrations of topological spaces.

\begin{defn}
    Let $X$ and $Y$ be sets, and let $R \subseteq X^{\ast} \times X^{\ast}$ and $S \subseteq Y^{\ast} \times Y^{\ast}$ be sets of relations. We call a map $f: X \to Y$ \emph{essential} with respect to the sets of relations $R$ and $S$ if

\begin{itemize}
    \nameditem{cond:E1}{(E1)} $f$ is surjective.
    \nameditem{cond:E2}{(E2)} $(v,w) \in R \Rightarrow (f^{\ast}(v),f^{\ast}(w)) \in S$ for all $v,w \in X^{\ast}$.
    \nameditem{cond:E3}{(E3)} If $(f^{\ast}(u),v') \in S$ holds for $u \in X^{\ast}$, $v' \in Y^{\ast}$, then there is a $v \in X^{\ast}$ with $(u,v) \in R$ and $f^{\ast}(v) = v'$.
\end{itemize}
\end{defn}

Let $M = \genrel{X}{R}_{\mathrm{mon}}$ and $N = \genrel{Y}{S}_{\mathrm{mon}}$. Note that by \labelcref{cond:E1} and \labelcref{cond:E2}, \cref{lem:respect_relations} implies that each essential map $f: X \to Y$ induces a well-defined surmorphism
\[
\overline{f}: M \twoheadrightarrow N; \ \overline{u} \mapsto \overline{f^{\ast}(u)}.
\]

Essentiality allows us to lift equivalences in the following sense:

\begin{pro} \label{pro:lifting_property_free_monoids}
    Given sets $X$ and $Y$, let $R$ and $S$ be sets of relations in $X^{\ast}$ and $Y^{\ast}$, respectively. Let $f: X \to Y$ be an essential map with respect to the sets $R$ and $S$.

    If $f^{\ast}(u) \approx_S v'$ for $u \in X^{\ast}$, $v' \in Y^{\ast}$, then there is a $v \in X^{\ast}$ with $u \approx_R v$ and $f^{\ast}(v) = v'$.
\end{pro}

\begin{proof}
    It is sufficient to show that $f^{\ast}(u) \sim_S v'$ for $u \in X^{\ast}$, $v' \in Y^{\ast}$, implies the existence of a $v \in X^{\ast}$ with $u \sim_R v$ and $f^{\ast}(v) = v'$.

    Since $f^{\ast}$ is obtained by substituting letters in $X$ by letters in $Y$, the relation $f^{\ast}(u) \sim_S v'$ implies that there are factorizations $u = u_1u_2u_3$, $v' = v_1^{\prime} v_2^{\prime} v_3^{\prime}$ such that $f^{\ast}(u_i) = v_i^{\prime}$ for $i \in \{1,3\}$ and such that $(f^{\ast}(u_2),v_2^{\prime}) \in S$. But by \labelcref{cond:E3}, there is a $v_2 \in X^{\ast}$ such that $f^{\ast}(v_2) = v_2^{\prime}$ and $(u_2,v_2) \in R$. Letting $v = u_1v_2u_3$, we have $u \sim_R v$ and $f^{\ast}(v) = f^{\ast}(u_1)f^{\ast}(v_2)f^{\ast}(u_3) = v_1^{\prime} v_2^{\prime} v_3^{\prime} = v'$.
\end{proof}

\begin{lem} \label{lem:lifting_property_monoids}
    Let $M = \genrel{X}{R}_{\mathrm{mon}}$ and $N = \genrel{Y}{S}_{\mathrm{mon}}$, and let $f: X \to Y$ be an essential map with respect to $R$ and $S$.
    
    Then for all $a',b'\in N$, $c \in M$ with $a'b' = \overline{f}(c)$ there are $a,b \in M$ with $\overline{f}(a) = a'$, $\overline{f}(b) = b'$ and $ab = c$.
\end{lem}

\begin{proof}
    Choose representatives $u',v' \in Y^{\ast}$ and $w \in X^{\ast}$ such that $\overline{u'} = a'$, $\overline{v'} = b'$ and $\overline{w} = c$. Then the equation $a'b' = \overline{f}(c)$ in $N$ is equivalent to the congruence $u'v' \approx_S f^{\ast}(w)$ in $Y^{\ast}$. By \cref{pro:lifting_property_free_monoids}
    this implies the existence of a word $z \in X^{\ast}$ with $w \approx_R z$ and $f^{\ast}(z) = u'v'$. But again, as $f^{\ast}$ acts by substituting letters in $X$ by letters in $Y$, we have a factorization $z = uv$ with $f^{\ast}(u) = u'$ and $f^{\ast}(v) = v'$.

    But with $a = \overline{u}$ and $b = \overline{v}$, we obtain $\overline{f}(a) = \overline{u'} = a'$ and, similarly, $\overline{f}(b) = b'$. Furthermore, since $uv \approx_R w$, we obtain that $ab = \overline{u} \overline{v} = \overline{w} = c$.
\end{proof}

Recall that for an element $x \in X^{\ast}$, its \emph{length} $\mathrm{len}(x) \in \N_0$ is defined by
\[
\mathrm{len}(x) = \begin{cases}
    0 & x = \varepsilon \\
    k & x = x_1x_2\ldots x_k \in X^k.
\end{cases}
\]
A relation $(v,w) \in X^{\ast} \times X^{\ast}$ is called \emph{homogeneous} if $\mathrm{len}(v) = \mathrm{len}(w)$, and it is easy to check that if $R \subseteq X^{\ast} \times X^{\ast}$ is a set of homogeneous relations, then $\overline{\mathrm{len}}(\overline{v}) = \mathrm{len}(v)$ defines a well-defined length function on $\genrel{X}{R}_{\mathrm{mon}}$. With this, we can now state:

\begin{thm} \label{thm:upp_under_essential_homomorphisms}
    Let $M = \genrel{X}{R}_{\mathrm{mon}}$ and $N = \genrel{Y}{S}_{\mathrm{mon}}$, where $R$ is a set of homogeneous relations, and let $f: X \to Y$ be an essential map with respect to $R$ and $S$. If $X$ is finite and $M$ satisfies the UPP, then $N$ satisfies the UPP.
\end{thm}

\begin{proof}
    By \labelcref{cond:E1}, $Y$ is finite. We first show that $S$ is a set of homogeneous relations. Given a relation $(v,w) \in S$, $\labelcref{cond:E1}$ allows us to write $v = f^{\ast}(v')$ for some $v' \in X^{\ast}$. But then there is a $w' \in X^{\ast}$ with $(v',w') \in R$ and $f^{\ast}(w') = w$. But $(v',w') \in R$ is homogeneous and $f^{\ast}$ preserves the length of words, thus $(v,w) = (f^{\ast}(v'),f^{\ast}(w'))$ is homogeneous as well.

    Suppose that $N$ does \emph{not} satisfy the UPP, and let the finite sets $\emptyset \neq C,D \subseteq N$ witness the failure of the UPP, meaning that for any $(c,d) \in C \times D$, there is an element $(c',d') \in C \times D$ with $(c,d) \neq (c',d')$ and $cd = c'd'$.

    $\overline{f}: M \to N$ is length-preserving, and there are only finitely many elements of given length in $M$ and $N$, thus it follows that $\overline{f}$ is finite-to-one. In particular, the sets $A = \overline{f}^{-1}(C)$ and $B = \overline{f}^{-1}(D)$ are finite. Furthermore, as $\overline{f}$ is surjective due to \labelcref{cond:E1}, we have $\overline{f}(A) = C$ and $\overline{f}(B) = D$. In particular, $A,B \neq \emptyset$.

    We will now show that $A,B$ witness that $M$ fails to satisfy the UPP. Let $(a,b) \in A \times B$. Since $(\overline{f}(a),\overline{f}(b)) \in C \times D$, there is $(c,d) \in C \times D$ with $(\overline{f}(a),\overline{f}(b)) \neq (c,d)$ and $\overline{f}(a)\overline{f}(b) = cd$. In particular, $\overline{f}(ab) = cd$. Now by \cref{lem:lifting_property_monoids}, there are $a',b' \in M$ such that $\overline{f}(a') = c$, $\overline{f}(b') = d$ and $a'b' = ab$. These satisfy $a' \in \overline{f}^{-1}(c) \subseteq \overline{f}^{-1}(C)$, $b' \in \overline{f}^{-1}(d) \subseteq \overline{f}^{-1}(D)$.

    Furthermore, as $(\overline{f}(a'),\overline{f}(b')) = (c,d) \neq (\overline{f}(a),\overline{f}(b))$, it follows that $(a,b) \neq (a',b')$. This proves that $M$ does not have the UPP.
\end{proof}

\section{Structure groups and the UPP} \label{sec:structure_groups_and_the_upp}

We now apply the methods from \cref{sec:essential_maps} to structure groups of set-theoretic solutions to the Yang--Baxter equation. This is done by exploiting the Garside structure on structure groups (\cref{thm:chouraqui_garside}).

\begin{pro} \label{pro:upp_in_garside_groups}
    Let $G$ be a Garside group with Garside monoid $G^+ \subseteq G$. Then $G^+$ has the UPP if and only if $G$ has the UPP.
\end{pro}

\begin{proof}
    It is readily checked that the UPP is inherited by submonoids. Therefore, we only have to show that $G$ satisfies the UPP if $G^+$ does.

    To this end, let $\emptyset \neq A,B \subseteq G$ be finite. By \cref{lem:shifting_subsets_into_G+}, there are integers $k,l$ such that $\Delta^kA, B\Delta^l \subseteq G^+$. Now the UPP for $G^+$ implies that there is a pair $(a,b) \in A \times B$ such that $(\Delta^ka)(b\Delta^l) \neq (\Delta^ka')(b'\Delta^l)$ for all $(a',b') \in A \times B$ with $(\Delta^ka,b\Delta^l) \neq (\Delta^ka',b'\Delta^l)$. But this is equivalent to $ab \neq a'b'$ for all $(a',b') \in A \times B$ with $(a,b) \neq (a',b')$. Since $A,B$ are arbitrary, this shows that $G$ satisfies the UPP.
\end{proof}

We can now prove:

\begin{pro} \label{pro:upp_for_factors_monoid}
    Let $(X,r)$, $(Y,s)$ be set-theoretic solutions to the Yang--Baxter equation, and let $f: X \twoheadrightarrow Y$ be a surmorphism of solutions. If $M(X,r)$ satisfies the UPP, then $M(Y,s)$ also satisfies the UPP.
\end{pro}

\begin{proof}
    If $(X,r)$ and $(Y,s)$ are solutions, we recall that $M(X,r) = \genrel{X}{R}_{\mathrm{mon}}$ where
    \[
    R = \{ (wx,yz) : w,x,y,z \in X;\ r(w,x) = (y,z) \}.
    \]
    Similarly, $M(Y,s) = \genrel{Y}{S}_{\mathrm{mon}}$ where
    \[
    S = \{ (w'x',y'z') : w',x',y',z' \in Y;\ s(w',x') = (y',z') \}.
    \]
    Note that the involutivity of $r$ and $s$ imply that $R$ and $S$ are symmetric.
    
    We now show that each surmorphism of solutions $f: (X,r) \twoheadrightarrow (Y,s)$ is also essential with respect to $R$ and $S$.
    
    By definition, $f$ is surjective; therefore, \labelcref{cond:E1} holds.
    
    Given $(wx,yz) \in R$, we have $r(w,x) = (y,z)$. Hence,
    \[
    s(f(w),f(x)) = (f \times f)(r(w,x)) = (f(y),f(z)),
    \]
    which implies that $(f^{\ast}(wx),f^{\ast}(yz)) \in S$. This confirms \labelcref{cond:E2}.

    In order to check \labelcref{cond:E3}, assume that $(f^{\ast}(wx), y'z') \in S$. Then
    \[
    (y',z') = s(f(w),f(x)) = (f \times f)(r(w,x)),
    \]
    and with $r(w,x) = (y,z)$, we have $(y',z') = (f(y),f(z))$. Therefore, $f^{\ast}(yz) = y'z'$ and $(wx,yz) \in R$.

    Therefore, we have shown that $f$ is essential. If $M(X,r)$ satisfies the UPP, then by \cref{thm:upp_under_essential_homomorphisms}, $M(Y,s)$ also satisfies the UPP.
\end{proof}

\begin{thm} \label{thm:upp_for_factors_group}
    Let $(X,r)$, $(Y,s)$ be set-theoretic solutions to the Yang--Baxter equation, and let $f: X \twoheadrightarrow Y$ be a surmorphism of solutions. If $G(X,r)$ satisfies the UPP, then $G(Y,s)$ also satisfies the UPP.
\end{thm}

\begin{proof}
    Suppose that $G(X,r)$ satisfies the UPP. $M(X,r)$ embeds as a Garside monoid into $G(X,r)$ due to \cref{thm:chouraqui_garside}. Therefore, the monoid $M(X,r)$ also satisfies the UPP and, by an application of \cref{pro:upp_for_factors_monoid}, so does $M(Y,s)$. Now \cref{pro:upp_in_garside_groups} proves that $G(Y,s)$ also satisfies the UPP.
\end{proof}

Recall that the \emph{Promislow group} is the group defined by generators and relations as
\[
P = \genrel{x,y}{x^2yx^2 = y;\ y^2xy^2 = x}_{\mathrm{gr}}.
\]

The Promislow group is a standard example of a torsionfree group without the UPP.

\begin{thm}
    $P$ is a torsionfree group without the UPP.
\end{thm}

\begin{proof}
    See \cite{Promislow_UPP}.
\end{proof}

Furthermore, we have the following result of Acri, Lutowski and Vendramin.

\begin{thm} \label{thm:ALV}
    Let $(X,r)$ be a solution to the Yang--Baxter equation of size $|X| \leq 8$ that is not of finite multipermutation level. Then $G(X,r)$ contains a subgroup isomorphic to $P$ if and only if $(X,r)$ is \emph{not} isomorphic to a solution of the form $(Y,s)$, where $s(x,y) = (\sigma_x(y),\tau_y(x))$ is from \cref{tab:structure_groups_without_promislow}\footnote{Table \ref{tab:structure_groups_without_promislow} has been reproduced from \cite{ALV_UPP} with the friendly permission from Leandro Vendramin.}.
\end{thm}

\begin{proof}
    \cite[Theorem 4.7]{ALV_UPP}.
\end{proof}

\begin{table}
\begin{center}
\begin{tabular}{|c|c|c||c|c|} 
\hline
$x$ & $\sigma_x$ & $\tau_x$ & $\sigma_x$ & $\tau_x$\\
\hline
1 & $(3745)$ & $(3648)$ & $(3745)(68)$ & $(3648)(57)$\\
2 & $(3745)$ & $(3648)$ & $(3745)(68)$ & $(3648)(57)$\\
3 & $(1826)$ & $(1527)$ & $(1826)(57)$ & $(1527)(68)$\\
4 & $(1826)$ & $(1527)$ & $(1826)(57)$ & $(1527)(68)$\\
5 & $(13872465)$ & $(16542873)$ & $(1465)(2387)$ & $(1654)(2873)$\\
6 & $(17842563)$ & $(13562478)$ & $(1784)(2563)$ & $(1478)(2356)$\\
7 & $(13872465)$ & $(16542873)$ & $(1465)(2387)$ & $(1654)(2873)$\\
8 & $(17842563)$ & $(13562478)$ & $(1784)(2563)$ & $(1478)(2356)$\\
\hline
1 & $(12)(4675)$ & $(12)(3685)$ & $(12)(35)(4867)$ & $(12)(3857)(46)$\\ 
2 & $(12)(4675)$ & $(12)(3685)$ & $(12)(35)(4867)$ & $(12)(3857)(46)$\\
3 & $(1435)(2786)$ & $(1578)(2643)$ & $(16582437)$ & $(17652843)$\\
4 & $(1587)(2634)$ & $(1345)(2876)$ & $(17562834)$ & $(15682347)$\\
5 & $(1823)(56)$ & $(1724)(56)$ & $(16582437)$ & $(17652843)$\\
6 & $(1823)(56)$ & $(1724)(56)$ & $(17562834)$ & $(15682347)$\\
7 & $(1587)(2634)$ & $(1345)(2876)$ & $(1325)(46)(78)$ & $(1426)(35)(78)$\\
8 & $(1435)(2786)$ & $(1578)(2643)$ & $(1325)(46)(78)$ & $(1426)(35)(78)$\\
\hline
1 & $(12)(78)$ & $(14)(67)$ & $(12)(35)(46)(78)$ & $(14)(28)(35)(67)$\\
2 & $(1584)(2673)$ & $(1265)(3784)$ & $(1324)(5867)$ & $(1243)(5786)$\\
3 & $(34)(56)$ & $(23)(58)$ & $(17)(28)(34)(56)$ & $(17)(23)(46)(58)$\\
4 & $(1485)(2376)$ & $(1562)(3487)$ & $(1423)(5768)$ & $(1342)(5687)$\\
5 & $(34)(56)$ & $(23)(58)$ & $(17)(28)(34)(56)$ & $(17)(23)(46)(58)$\\
6 & $(1485)(2376)$ & $(1562)(3487)$ & $(1423)(5768)$ & $(1342)(5687)$\\ 
7 & $(12)(78)$ & $(14)(67)$ & $(12)(35)(46)(78)$ & $(14)(28)(35)(67)$\\
8 & $(1584)(2673)$ & $(1265)(3784)$ & $(1324)(5867)$ & $(1243)(5786)$\\
\hline
1 & $(13687542)$ & $(13867524)$& $(1652)$ & $(1854)$\\
2 & $(17)(2583)(46)$ & $(1278)(35)(46)$ & $(17645328)$ & $(12835647)$\\
3 & $(18657243)$ & $(16857423)$ & $(3874)$ & $(2367)$\\
4 & $(1476)(28)(35)$ & $(17)(28)(3456)$ & $(14635827)$ & $(17825346)$\\
5 & $(18657243)$ & $(16857423)$ & $(1652)$ & $(1854)$\\
6 & $(1476)(28)(35)$ & $(17)(28)(3456)$ & $(17645328)$ & $(12835647)$\\
7 & $(13687542)$ & $(13867524)$ & $(3874)$ & $(2367)$\\
8 & $(17)(2583)(46)$ & $(1278)(35)(46)$ & $(14635827)$ & $(17825346)$\\
\hline
\end{tabular}
\caption{Representatives of all isomorphism classes of solutions $(Y,s)$ with $|Y| \leq 8$, that are not of finite multipermutation level, where $G(Y,s)$ does not contain the Promislow group.} \label{tab:structure_groups_without_promislow}
\end{center}
\end{table}

We can now prove the following theorem.

\begin{thm}
    Let $(X,r)$ be a set-theoretic solution to the Yang--Baxter equation with $|X| \leq 8$. Then $(X,r)$ is of finite multipermutation level if and only if $G(X,r)$ satisfies the UPP.

    In particular, there exist set-theoretic solutions $(X,r)$ such that $G(X,r)$ does not satisfy the UPP, while not containing an isomorphic copy of the Promislow group.
\end{thm}

\begin{proof}
    If $(X,r)$ is of finite multipermutation level, then $G(X,r)$ is left-orderable by \cite[Theorem 2.1]{BCV_Left_Orderable}, so $G(X,r)$ also satisfies the UPP.
    
    Now suppose that $(X,r)$ is not of finite multipermutation level. By \cref{thm:ALV}, if $(X,r)$ is not isomorphic to one of the solutions in \cref{tab:structure_groups_without_promislow}, $G(X,r)$ contains a subgroup isomorphic to the Promislow group and therefore does not satisfy the UPP. On the other hand, if $(X,r)$ is isomorphic to one of the solutions from \cref{tab:structure_groups_without_promislow}, then it is quickly checked that $(X,r)$ is retractable. In that case, $(X^{\prime},r^{\prime}) = \Ret(X,r)$ is of size $|X^{\prime}| < 8$, and is not of finite multipermutation level. Therefore, by \cref{thm:ALV}, $G(X^{\prime},r^{\prime})$ contains a copy of the Promislow group and does not satisfy the UPP.

    As the canonical map $X \twoheadrightarrow X'$; $x \mapsto [x]$ provides a surmorphism from $(X,r)$ to $\Ret(X,r)$ (see \cref{subsec:the_YBE}), it follows by an application of \cref{thm:upp_for_factors_group} that $G(X,r)$ does not satisfy the UPP.

    Furthermore, if $(X,r)$ is from \cref{tab:structure_groups_without_promislow}, then $G(X,r)$ does not contain an isomorphic copy of $P$. Therefore, each such solution provides a structure group without the UPP that does not contain an isomorphic copy of $P$ as a subgroup.
\end{proof}

Finally, we show that the UPP for $G(X,r)$ is also inherited by subfactors of solutions.

Given two solutions $(X,r)$ and $(Y,s)$, we write $(X,r) \preceq (Y,s)$ if and only if there is a subsolution $(Z,s) \subseteq (Y,s)$, together with a surmorphism of solutions $f: (Z,s) \twoheadrightarrow (X,r)$.

\begin{pro}
    The relation $\preceq$ is a partial order relation on the set of isomorphism classes of finite solutions.
\end{pro}

\begin{proof}
    The relation $\preceq$ is reflexive, because $(X,r) \subseteq (X,r)$ is a subsolution and $\id_X: X \twoheadrightarrow X$ is a surmorphism.

    It is also antisymmetric: suppose that $(X,r) \preceq (Y,s)$ and $(Y,s) \preceq (X,r)$. As $(X,r)$ is a factor of a subsolution of $(Y,s)$, we have $|X| \leq |Y|$. Similarly, $|Y| \leq |X|$; therefore $|X| = |Y|$. Now if $(Z,s) \subseteq (Y,s)$ is a subsolution and $f: (Z,s) \twoheadrightarrow (X,r)$ is a surmorphism, this implies that $Z = Y$ and that $f$ is an isomorphism.

    We now show that $\preceq$ is also transitive. Suppose that $(X,r) \preceq (Y,s) \preceq(Z,t)$. Then there is a subsolution $(A,t) \subseteq (Z,t)$ together with a surmorphism $f: (A,t) \twoheadrightarrow (Y,s)$, and a subsolution $(B,s) \subseteq (Y,s)$ together with a surmorphism $g: (B,s) \twoheadrightarrow (X,r)$. It is easily checked that $(f^{-1}(B),t) \subseteq (Z,t)$ is a subsolution and $gf|_{f^{-1}(B)}: (f^{-1}(B),t) \twoheadrightarrow (X,r)$ is a surmorphism.
\end{proof}

We can now present a slight generalization of \cref{thm:upp_for_factors_group}.

\begin{thm}
    Suppose that $(X,r)$, $(Y,s)$ are solutions to the Yang--Baxter equation with $(X,r) \preceq (Y,s)$. If $G(Y,s)$ satisfies the UPP, then $G(X,r)$ satisfies the UPP.
\end{thm}

\begin{proof}
    Let $(X,r)$ and $(Y,s)$ be as in the statement of the theorem, and suppose that $(Z,s) \subseteq (Y,s)$ is a subsolution such that there is a surmorphism $f: (Z,s) \twoheadrightarrow (X,r)$. By \cref{thm:upp_for_factors_group}, it is sufficient to show that $G(Z,s)$ satisfies the UPP. By the same argument as in the proof of \cref{thm:upp_for_factors_group}, it is sufficient to show that $M(Z,s)$ satisfies the UPP. As $M(Y,s)$ satisfies the UPP, this will follow from showing that $M(Z,s) \cong \langle Z \rangle$ where $\langle Z \rangle$ denotes the submonoid generated by $Z$ in $M(Y,s)$. Since $s(Z^2) = Z^2$, it follows that two expressions in $Z^{\ast}$ are equivalent with respect to the defining relations of $M(Z,s)$ if and only if they are equivalent with respect to the defining relations of $M(Y,s)$.
\end{proof}

Note that our method yields examples virtually abelian groups without the UPP that are of the form $G(X,r)$ and do not contain the Promislow group as a subgroup. This motivates the following question.

\begin{problem}
    Let $k$ be a field, and let $G = G(X,r)$ be one of the groups constructed in this work. Does the group ring $k[G]$ contain nontrivial units? If so, describe them explicitly!
\end{problem}

\bibliography{references}
\bibliographystyle{abbrv}

\end{document}